\magnification=\magstep1
\input amstex
\documentstyle {amsppt}
\NoBlackBoxes
\NoRunningHeads
\nologo
\pageheight {8.5 truein}
\pagewidth {6.5 truein}

\topmatter

\title
On a class of Koszul algebras associated to directed graphs 
\endtitle
\author
Vladimir Retakh \footnote {partially supported by NSA \hfill \hfill}, Shirlei Serconek 
 and Robert Lee Wilson
\endauthor
\address
\newline
V.R., R.W.: Department of Mathematics, Rutgers University,
Piscataway,
NJ 08854-8019
\newline
S.S:
IME-UFG
CX Postal 131
Goiania - GO
CEP 74001-970 Brazil
\endaddress
\email
\newline vretakh\@math.rutgers.edu, serconek\@math.rutgers.edu, rwilson\@math.rutgers.edu
\endemail

\endtopmatter

\endtopmatter
\document

\centerline {\bf Abstract}
\smallskip 
In math.QA/0506507 I. Gelfand and the authors of this paper introduced
a new class of algebras associated to directed graphs. In this
paper we show that these algebras are Koszul for a large class
of layered graphs.

\bigskip
 
\head \bf 0. Introduction  \endhead

In \cite {GRSW} I. Gelfand and the authors of this paper associated
to any layered graph $\Gamma $ an algebra $A(\Gamma )$ and
constructed a basis in $A(\Gamma )$ when the graph is a layered graph
with a unique minimal vertex. 

The algebra $A(\Gamma )$ is a natural generalization of
universal algebra $Q_n$ of pseudo-roots of noncommutative polynomials
introduced in \cite {GRW}.
In fact, $A(\Gamma )$ is isomorphic to $Q_n$ when
$\Gamma $ is the hypercube of dimension $n$, i.e. the graph of all subsets
of a set with $n$ elements.

The algebras $Q_n$  have a rich and interesting structure
related to factorizations of polynomials over noncommutative rings. On one
hand,
$Q_n$ is a ``big algebra" (in particular, it contains free
subalgebras on several generators and so  has
an exponential growth). On the
other hand, it is rather ``tame": it is a quadratic algebra, one can
construct a linear
basis in $Q_n$ \cite {GRW},  compute its Hilbert series \cite {GGRWS}, prove that $Q_n$ is Koszul  \cite {SW, Pi}, and construct interesting quotients of $Q_n$ \cite {GGR}.

Since the algebra $A(\Gamma )$ is a natural generalization of $Q_n$
one would expect that for a ``natural" class of graphs the algebra $A(\Gamma
)$
is Koszul. In this paper we prove  this assertion when $\Gamma$ is a{ \it uniform layered graph}; see
Definition 3.3. The Hasse graph of ranked modular lattices with a unique minimal element is an example of such a graph.

Compared to the proof given in \cite {SW} for the algebra $Q_n$, our proof is
much simpler, and
more geometric.

\head \bf 1. The algebra $A(\Gamma)$ as a quotient of $T(V^+)$
\endhead
We begin by recalling (from \cite{GRSW}) the definition of the algebra $A(\Gamma).$
Let  $\Gamma = (V, E)$ be a {\bf directed graph}. That is, $V$ is a set (of vertices), 
$E$ is a set (of edges), and ${\bold t}: E \rightarrow V$ and ${\bold h}: E \rightarrow V$ are functions. (${\bold t}(e)$ is the {\it tail} 
of $e$ and ${\bold h}(e)$ is the {\it head} of $e$.)       
\smallskip

We say that $\Gamma$ is {\bf layered} if $V = \cup _{i=0}^n V_i$, $E = \cup_{i=1}^n E_i$, 
 ${\bold t}: E_i \rightarrow V_i$, \ ${\bold h}: E_i \rightarrow V_{i-1}$. Let $V^+ = \cup_{i=1}^n V_i.$
\smallskip  
We will assume throughout the remainder of the paper that $\Gamma = (V, E)$ is a layered graph with 
$V = \cup_{i=0}^n V_i$, that $V_0 = \{*\}$, and that, for every $v \in  V^+$, 
$\{ e \in E \ | \ {\bold t}(e) = v \} \neq \emptyset$. 

If $v, \ w \in V$, a {\bf path} from $v$ to $w$ is a sequence of edges 
$\pi = \{ e_1, e_2, ...,e_m \}$ with ${\bold t}(e_{1}) = v$, ${\bold h}(e_m) = w$ and ${\bold t}(e_{i+1}) = {\bold h}(e_i)$ for 
$1 \leq i < m$.  We write $v = {\bold t}(\pi)$, $w = {\bold h}(\pi)$.  We also write $v > w$ if there  is a path from 
$v$ to $w$.  Define $$P_{\pi}(t) = (1 -te_1)(1-te_2)...(1-te_m) \in T(E)[t]/(t^{n+1})$$ and write
$$P_{\pi}(t) = \sum_{j=0}^n e(\pi,j)t^j.$$

Recall (from \cite{GRSW}) that $R$ denotes the ideal of $T(E)$, the tensor algebra on $E$ over the field $F$, generated by $$\{ e(\pi_1,k) - e(\pi_2,k) \ | \ {\bold t}(\pi_1)={\bold t}(\pi_2), {\bold h}(\pi_1)={\bold h}(\pi_2) , \ 1 \leq k \leq l(\pi_1) \}.$$
Also, by Lemma 2.5 of \cite{GRSW}, $R$ is actually generated by the smaller set 
$$\{ e(\pi_1,k) - e(\pi_2,k) \ | \ {\bold t}(\pi_1)={\bold t}(\pi_2), {\bold h}(\pi_1)={\bold h}(\pi_2) = *, \ 1 \leq k \leq l(\pi_1) \}.$$

\definition{Definition 1.1} $A(\Gamma) = T(E)/R.$
\enddefinition

Note that $e(\pi,k) \in T(E)_k$.  Thus $R = \sum_{j=1}^\infty R_j$ is a graded ideal in $T(E).$ 
We write $\tilde {e}(\pi,k)$ for the image of $e(\pi,k)$ in $A(\Gamma)$.

In fact, $A(\Gamma)$ may also be expressed as a quotient of $T(V^+)$.  To verify this 
we  need a general result about quotients of the tensor algebra by graded ideals.  
Let $W$ be a vector space over $F$ and let $I = \sum_{j=1}^\infty I_j$ be a graded
ideal in the tensor algebra $T(W)$.  Let $\psi$ denote the canonical map from $W$ to the quotient space
$W/I_1$ and let $<I_1>$ denote the ideal in $T(W)$ generated by $I_1$. Then 
$\psi$ induces  a surjective homomorphism of graded algebras
$$\phi:T(W) \rightarrow T(W/I_1) \cong T(W)/<I_1>.$$  Consequently, by the Third Isomorphism Theorem, we have:

\proclaim {Proposition 1.2} 
$$T(W)/I \cong T(W/I_1)/\phi(I)$$
where $\phi(I)  =\sum_{j=2}^\infty \phi(I_j)$ is a graded ideal of $T(W/I_1).$
\endproclaim

We now apply this to the presentation of the algebra $A(\Gamma).$  
Recall that for each vertex $v \in V^+$ 
there is a distinguished edge $e_v$ with ${\bold t}(e_v) = v.$ 
Recall further that  for $v \in V^+$ 
we define $v^{(0)} = v$ and 
$v^{(i+1)} = {\bold h}(e_{v^{(i)}})$ for $0 \le i <|v|-1$  
and that we set 
$e(v,1) = e_{v^{(0)}} + e_{v^{(1)}} + ... + e_{v^{(|v|-1)}}.$  
Thus 
$$e_v = e(v,1) - e(v^{(1)},1) = e({\bold t}(e_v),1) - e({\bold h}(e_v),1).$$
Let $E' = \{e_v |v \in V^+\}.$ 
Define $\tau:FE \rightarrow FE'$ by 
$$\tau(f) = e({\bold t}(f),1) - e({\bold h}(t),1).$$
Then $\tau$ is a projection of $FE$ onto $FE'$ with kernel $R_1$.  

Now define $\eta: FE' \rightarrow FV^+$ by
$$\eta: e_v \mapsto v.$$
Then $\eta$ is an isomorphism of vector spaces and $\eta\tau$ induces an isomorphism 
$$\nu: FE/R_1 \rightarrow FV^+.$$
As above, $\nu$ induces a surjective homomorphism of graded algebras
$$\theta: T(E) \rightarrow T(V^+).$$

Thus Proposition 1.2 gives:
\proclaim {Corollary 1.3} $A(\Gamma) \cong T(V^+)/\theta(R).$
\endproclaim

It is important to write generators for the ideal $\theta(R)$ explicitly.  
Since $R$ is generated by $R_1$ together with the 
elements of the form $e(\pi_1,k) - e(\pi_2,k)$ it will be sufficient to write $\theta(e(\pi,k))$ explicitly.
Let $\pi = \{e_1,e_2,...,e_m\}$ be a path, let ${\bold{t}}(e_i) = v_{i-1}$ for $1 \le i \le m$ and let ${\bold h}(e_m ) = v_m.$  Then
$$e(\pi,k) = (-1)^k\sum_{1 \le i_1 < ... < i_k \le m} e_{i_1}...e_{i_k}.$$  
Now 
$\nu(e_i) = e_{v_{i-1}} - e_{v_i}$  and so 
$\eta\nu(e_i) = v_{i-1} - v_i.$  
Since $\theta$ is induced by $\eta\nu$  we have:
\proclaim{Lemma 1.4}
$$\theta(e(\pi,k)) = (-1)^k\sum_{1 \le i_1 < ... < i_k \le m} (v_{i_1 - 1}- v_{i_1})...(v_{i_k - 1 }-v_{i_k}).$$
\endproclaim

\bigskip
\head \bf 2. A presentation of $ gr \ A(\Gamma)$
\endhead
Let $W=\sum_{k=0}^\infty W_k$ be  a graded vector space. We begin by recalling some basic properties of $T(W).$

$T(W)$ is bi-graded,  that is, in addition to the usual grading (by degree in the tensor algebra), there is another grading induced by the grading of $W$.  Thus
$$T(W) = \sum_{i=0}^\infty T(W)_{[i]}$$
where
$$T(W)_{[i]} = span\{w_1...w_r|r \ge 0, w_j \in W_{[l_j]}, l_1 + ... + l_r = i\}.$$
This grading induces a filtration on $T(W).$  Namely
$$T(W)_i = T(W)_{[i]} + T_{[i-1]} + ... + T(W)_{[0]} = $$
$$span \{ w_1...w_r \ | r \ge 0, w_j \in W_{[l_j]}, l_1 + ... + l_r \le i \}.$$
Since $T(W)_i/T(W)_{i-1} \cong T(W)_{[i]}$ 
we may identify $T(W)$ with its associated graded algebra.
Define  a map
$$gr:T(W) \rightarrow T(W) = gr \ T(W)$$
by
$$gr: \lambda \mapsto \lambda$$
for $\lambda \in F.1$ and
$$gr:u = \sum_{i=0}^k u_k \mapsto u_k$$
where $k > 0, u_i \in T(W)_{[i]}$ and $u_k \ne 0.$

\proclaim{Lemma 2.1} Let $W$ be a graded vector space and $I$ be an ideal in $T(W).$  Then
$$gr \ (T(W)/I) \cong T(W)/(gr \ I).$$
\endproclaim
\demo {Proof} We have $$(gr\ I)_{[k]}=T(W)_{[k]}\cap (T(W)_{k-1} +I).$$
Therefore
$$gr\ (T(W)/I)_{[k]}=
(T(W)/I)_k/(T(W)/I)_{k-1}=
((T(W)_k+I)/I)/((T(W)_{k-1}+I)/I)$$
$$\cong
(T(W)_k+I)/(T(W)_{k-1}+I)=
(T(W)_{[k]}+T(W)_{k-1}+I)/(T(W)_{k-1}+I)$$
$$\cong T(W)_{[k]}/(T(W)_{[k]}\cap (T(W)_{k-1}+I))
=T(W)_{[k]}/(gr\ I)_{[k]}.$$
\hfill $\square$ \enddemo

The decomposition of $V$ into layers induces a grading of the vector space $FV^+$.  Thus the tensor algebra $T(V^+)$ is graded and filtered as above. The following lemma shows that 
this filtration on $T(V^+)$ agrees with that induced by the filtration on $T(E)$.

\proclaim{Lemma 2.2} For all $i \ge 0, T(V^+)_i = \theta(T(E)_i).$
\endproclaim
\demo{Proof} This holds for $i = 0$ since $T(V^+)_0 = T(E)_0 = F.$  Furthermore,
$T(E)_1$ is spanned by $1$ and $\{f|f \in E_1\}$.  For $f \in E_1$ 
we have $\tau(f) = e({\bold t}(f), 1) - e({\bold h}(f),1)$, but ${\bold h}(f) = * $ and $e(*,1) = 0$ so
$\tau(f) = e({\bold t}(f),1).$  Hence $\eta\tau(f) = {\bold t}(f).$  Thus $T(V^+)_1 = \theta(T(E)_1).$

Now assume $T(V^+)_{i-1} = \theta(T(E)_{i-1})$.  Then $\theta(T(E)_i)$ is spanned by 
$$\theta(\{ e_1...e_r \ | \ r \geq 0, |e_1| + ... + |e_r| \leq i \}) = $$
$$\{ ({\bold t}(e_1) - {\bold h}(e_1))...({\bold t}(e_r) - {\bold h}(e_r)) \ | \ r \geq 0, |e_1| + ... + |e_r| \leq i \}.$$
Let $u = ({\bold t}(e_1) - {\bold h}(e_1))...({\bold t}(e_r) - {\bold h}(e_r)).$  
Then if $|e_1| + ... + |e_r| \le i$ 
we have $$u \equiv {\bold t}(e_1) ... {\bold t}(e_r) \ mod \ T(V^+)_{i-1}.$$
The lemma then follows by induction.
\hfill $\square$ \enddemo

\proclaim {Corollary 2.3} $A(\Gamma) \cong T(V^+)/\theta(R)$ as filtered algebras.
\endproclaim

If $u \in A(\Gamma)_i, u \notin A(\Gamma)_{i-1}$ we write $|u| = i.$

As before, let $\pi = \{e_1,e_2,...,e_m\}$ be a path and let ${\bold{t}}(e_i) = v_{i-1}$ for $1 \le i \le m$ and ${\bold h}(e_m ) = v_m.$  For $1 \le k \le m+1$ set
$$v(\pi,k) = v_0...v_{k-1}.$$

\proclaim{Lemma 2.4} Let $\pi_1,\pi_2$ be paths with ${\bold t}(\pi_1) = {\bold t}(\pi_2)$ and let 
$1 \le k \le l(\pi_1).$  Then
$$v(\pi_1,k) - v(\pi_2,k) \in gr \ \theta(R).$$
\endproclaim
\demo{Proof}  We may extend $\pi_1, \pi_2$ to paths $\mu_1, \mu_2$ such that ${\bold h}(\mu_1 ) = {\bold h}(\mu_2) = *.$
Then $e(\mu_1,k) - e(\mu_2,k) \in R.$  The result now follows from Lemma 1.4.
\hfill $\square$ \enddemo

Let $R_{gr}$ denote the ideal generated by 
$$\{v(\pi_1,k) - v(\pi_2,k)|{\bold t}(\pi_1) = {\bold t}(\pi_2), 2 \le k \le l(\pi_1)\}.$$

\proclaim{Proposition 2.5}  $gr \ A(\Gamma) \cong T(V^+)/R_{gr}.$
\endproclaim
\demo{Proof}  We begin by recalling the description of a basis for $gr \ A(\Gamma).$

We say that a pair $(v,k)$, $v \in V$, $0 \leq k \leq |v|$ can be {\it composed} with the pair $(u,l)$, 
$u \in V$, $0 \leq l \leq |u|$, if $v > u$ and $|u| = |v| - k$. If $(v,k)$ can be composed with $(u,l)$ 
we write $ (v,k) \models (u,l)$.
Let ${\bold B_1}(\Gamma)$   be the set of all sequences
$$ {\bold b} = ((b_1,m_1),(b_2,m_2),...,(b_k,m_k))$$  where 
$k \geq 0$, $b_1, b_2, ..., b_k \in V$, $0 \leq m_i \leq |b_i|$ for $1 \leq i \leq k$.
Let $${\bold B}(\Gamma) = 
 \{  {\bold b} = ((b_1,m_1),(b_2,m_2),...,(b_k,m_k)) \in {\bold B_1}(\Gamma) \ | \ \atop 
        (b_i,m_i) \not\models  (b_{i+1},m_{i+1}) ,   \  \  1  \leq i < k \}.$$

For  $${\bold b} = ((b_1,m_1),(b_2,m_2),...,(b_k,m_k)) \in {\bold B_1}(\Gamma)$$ 
set 
$$\tilde{e}({\bold b}) = \tilde{e}(b_1,m_1)...\tilde{e}(b_k,m_k).$$

Clearly
$\{ \tilde{e}({\bold b}) \ | \ {\bold b} \in {\bold B_1}(\Gamma) \}$ spans $A(\Gamma)$.
Writing
$\bar {e}({\bold b}) = \tilde {e}({\bold b}) + A(\Gamma)_{i-1} \in gr \ A(\Gamma)$
where $|\tilde {e}({\bold b})| = i$,   
Corollary 4.4 of \cite {GRSW} shows that 
$\{\bar {e}({\bold b})|{\bold b} \in {\bold B}(\Gamma)\}$ 
is a basis for $gr \ A(\Gamma)$.

Let $$\check{}: T(V^+) \rightarrow T(V^+)/R_{gr}$$
denote the canonical mapping. 
Write $\check{e}(b,m)$ for the image of $e(b,m)$ 
and $\check{e}({\bold b})$ for the image of $e({\bold b})$.

By Lemma 2.1 and Corolooary 1.3, we have $gr \ A(\Gamma) \cong T(V^+)/(gr \ \theta(R)).$ 
Since $R_{gr} \subseteq  \ gr \ \theta(R)$ (by Lemma 2.4) the canonical map
$$T(V^+) \rightarrow T(V^+)/(gr \ \theta(R))$$
induces a homomorphism 
$$\alpha: T(V^+)/R_{gr} \rightarrow gr \ A(\Gamma).$$
Clearly 
$$\alpha: \check{e}({\bold b}) \mapsto  \bar{e}({\bold b}).$$
Also, since $\{\bar {e}({\bold b})|{\bold b} \in {\bold B}(\Gamma)\}$ 
is a basis for $gr \ A(\Gamma)$, there is a linear map 
$$\beta: gr \ A(\Gamma) \rightarrow T(V^+)/R_{gr}$$
defined by
$$\beta: \bar{e}({\bold b}) \mapsto  \check{e}({\bold b}).$$
As $\alpha$ and $\beta$ are inverse mappings, the proof is complete.
\hfill $\square$ \enddemo

\bigskip
\head \bf 3. The quadratic algebra $A(\Gamma)$
\endhead

 We will now see that, for certain graphs $\Gamma$, $R$ is generated by $R_1 +R_2.$

\definition{Definition 3.1} Let $\Gamma$ be a layered graph and $v \in V_j, j \ge 2.$  For $1 \le i \le j$ define ${\Cal S}_i(v) = \{w \in V_{j-i}|v > w\}.$  
\enddefinition

\definition{ Definition 3.2} For $v \in V_j, j \ge 2$, let 
$\sim_v$ denote the equivalence relation on ${\Cal S}_1(v)$ generated by 
$u \sim_v w $ if ${\Cal S}_1(u) \cap {\Cal S}_1(w) \ne \emptyset.$
\enddefinition

\definition{ Definition 3.3} The layered graph $V$ is said to be {\bf uniform} if, for every $v \in V_j, j \ge 2$, all elements of ${\Cal S}_1(v)$ are equivalent under $\sim_v$.
\enddefinition

\proclaim {Lemma 3.4}  Let $\Gamma$ be a uniform layered graph.  Then $R$ is generated by $R_1 + R_2$, in fact, $R$ is generated by 
$R_1 \cup \{e(\tau_1,2) - e(\tau_2,2)|{\bold t}(\tau_1)= {\bold t}(\tau_2),  {\bold h}(\tau_1)= {\bold h}(\tau_2),
|\tau_1| = 2 \}.$

\endproclaim
\smallskip
\demo{Proof} Let $S$ denote the ideal of $T(E)$ generated by $R_1 \cup \{e(\tau_1,2) - e(\tau_2,2)|{\bold t}(\tau_1)= {\bold t}(\tau_2),  {\bold h}(\tau_1)= {\bold h}(\tau_2),
|\tau_1| = 2 \}.$

We must show that if $\pi_1, \pi_2$ are paths in 
$\Gamma$ with ${\bold t}(\pi_1) = {\bold t}(\pi_2),$ and
${\bold h}(\pi_1) = {\bold h}(\pi_2) = *$, 
then $P_{\pi_1}(t) - P_{\pi_2}(t) \in S[t]$, or, equivalently,
$P_{\pi_1}(t) \in (1 + S[t])P_{\pi_2}(t).$  
This is clear if $l(\pi_1) \le 2.$  We will 
proceed by induction on $l(\pi_1)$.  Thus we will assume 
that $k \ge 3$, that $l(\pi_1) = k,$ and that whenever $\mu_1, \mu_2$ are 
paths in $\Gamma$ with ${\bold t}(\mu_1) = {\bold t}(\mu_2),$ and
${\bold h}(\mu_1) = {\bold h}(\mu_2) = *$, and $l(\mu_1) < k$, then $P_{\mu_1}(t) - P_{\mu_2}(t) \in S[t].$

Write $\pi_1 = (e_1,e_2,..., e_k), \pi_2 = (f_1, f_2, ..., f_k)$.  
We first consider the special case in which ${\bold h}(e_1) > {\bold h}(f_2)$ 
(and so there is an edge, say $g$, with ${\bold t}(g) = {\bold h}(e_1), {\bold h}(g) = {\bold h}(f_2)).$  
Consequently, $P_{(e_1,g)}(t) \in (1 + S[t])P_{(f_1,f_2)}(t).$ Write $\pi_1 = (e_1,e_2)\nu_1$ 
and $\pi_2 = (f_1,f_2)\nu_2.$  Then
$$P_{\pi_1}(t) = (1 - te_1)(1 - te_2)P_{\nu_1}(t)$$
and
$$P_{\pi_2}(t) = (1 - tf_1)(1 - tf_2)P_{\nu_2}(t)$$
so
$$P_{\pi_2}(t) = $$
$$(1-tf_1)(1-tf_2)((1-tg)^{-1}(1-te_1)^{-1}(1-te_1)(1-tg))P_{\nu_2}(t) \times $$
$$P_{\nu_1}(t)^{-1}((1-te_2)^{-1}(1-te_1)^{-1}P_{\pi_1}(t)).$$
Let $\mu_1 = e_2\nu_1$ and $\mu_2 = g\nu_2$.  
Then, by the induction assumption, 
$$(1-tg)P_{\nu_2}(t)P_{\nu_1}(t)^{-1}(1-te_2)^{-1} = P_{\mu_2}(t)P_{\mu_1}(t)^{-1} \in 1 + S[t].$$  
Consequently, 
$$(1-te_1)(1-tg)P_{\nu_2}(t)P_{\nu_1}(t)^{-1}(1-te_2)^{-1}(1-te_1)^{-1}$$
$$ \in (1-te_1)(1 + S[t])(1 - te_1)^{-1} = 1 + S[t]$$ 
and so we have
$$P_{\pi_2}(t) \in (1 + S[t])P_{\pi_1}(t).$$

In the general case, let ${\bold h}(e_1) = u$ and ${\bold h}(f_1) = w$.  Then $u,w \in {\Cal S}_1(v)$ so, since $\Gamma$ is uniform, there exist $a_1,...,a_{l+1} \in {\Cal S}_1(v)$ with $a_1 = u, a_{l+1} = w$ and $ b_1,...,b_l \in V$ with
$b_i \in {\Cal S}_1(a_i) \cap {\Cal S}_1(a_{i+1})$ 
for $1 \le i \le l$.  For $1 \le i \le l$, let $\tau_i$ 
be a path from $b_i$ to *.  
For $2 \le i \le  l$ let $g_i \in E$ satisfy ${\bold t}(g_i) = {\bold t}(\pi_1), {\bold h}(g_i) = a_i.$
For $1 \le i \le l$ let $r_i \in E$ 
satisfy ${\bold t}(r_i) = a_i$ and ${\bold h}(r_i) = b_i$ and let $s_i \in E$ 
satisfy ${\bold t}(s_i) = a_{i+1}, {\bold h}(s_i )= b_i.$  Then the 
previously considered case shows that 
$$P_{\pi_1}(t) \in (1 + S[t])P_{g_2s_1\tau_1}(t);$$ 
$$P_{g_is_{i-1}\tau_{i-1}}(t) 
\in (1 + S[t]) P_{g_{i+1}s_i\tau_i}(t)$$ 
for $2 \le i \le l-1$;  
$$P_{g_ls_{l-1}\tau_{l-1}}(t) \in (1 + S[t])P_{f_1s_l\tau_l}(t);$$ 
and 
$$P_{f_1s_l\tau_l}(t) \in (1 + S[t])P_{\pi_2}(t),$$ 
proving the lemma.
\hfill $\square$ \enddemo

Now assume that $\Gamma$ is a uniform layered graph.  Then $R$ is generated
by $R_1 + R_2$, in fact, $R$ is generated by 
$R_1 \cup {\Cal R}_2$ where ${\Cal R}_2 = \{e(\tau_1,2) - e(\tau_2,2)|{\bold t}(\tau_1)= {\bold t}(\tau_2),  {\bold h}(\tau_1)= {\bold h}(\tau_2),
|\tau_1| = 2 \}.$
Set $R_V = <\theta({\Cal R}_2)>$.
\proclaim{Proposition 3.5} Let $\Gamma$ be a uniform layered graph.  Then $A(\Gamma) \cong T(V^+)/R_V$ is a quadratic algebra and $R_V$ is generated by $$\{v(u-w) - u^2 + w^2 + (u-w)x| v \in \cup_{i=2}^n V_i, u,w \in {\Cal S}_1(v), x \in {\Cal S}_1(u) \cap {\Cal S}_1(w)\}.$$
\endproclaim
\smallskip
\demo{Proof} By Lemma 3.4, $R_V$ is generated by 
$$\theta(\{e(\tau_1,2) - e(\tau_2,2)|{\bold t}(\tau_1)= {\bold t}(\tau_2),  {\bold h}(\tau_1)= {\bold h}(\tau_2),
|\tau_1| = 2 \}.$$
Let $\tau_1 = (e,f), \tau_2 = (e',f'), {\bold t}(e) = {\bold t}(e') = v, {\bold h}(e ) = u, {\bold h}(e') = w, {\bold h}(f)=
{\bold h}(f') = x.$  Then 
$$\theta(e(\tau_1,2) - e(\tau_2,2))= (v-u)(u-x) - (v-w)(w-x) = v(u-w) - u^2 + w^2 + (u-w)x.$$
\hfill $\square$ \enddemo

Combining this proposition with the results of the previous section 
we obtain the following presentation for $gr \ A(\Gamma))$.

\proclaim{Proposition 3.6} Let $\Gamma$ be a uniform layered graph.  
Then $gr \ A(\Gamma) \cong T(V^+)/R_{gr}$ is a quadratic algebra 
and $R_{gr}$ is generated by 
$$\{v(u-w) | v \in \cup_{i=2}^n V_i, u,w \in {\Cal S}_1(v), x \in {\Cal S}_1(u) \cap {\Cal S}_1(w)\}.$$
\endproclaim
\bigskip

\head \bf 4. $gr \ A(\Gamma)$ is a Koszul algebra
\endhead
\bigskip

If $W$ is a graded subspace of $V^2$ we write 
$$W^{(k)} = \bigcap_{i=0}^{k-2} V^iWV^{k-i-2}$$
so that
$$(gr \ W)^{(k)} = \bigcap_{i=0}^{k-2} V^i(gr \ W)V^{k-i-2} .$$
Then, by Proposition 3.6,
$$(gr \ R_{gr})^{(k)} \subset  span \{ v(\pi,k) \ | \  \pi \ is \ a \ path, \ l(\pi) \ge k \} $$ 

To simplify notation, we will write ${\bold V}$ for $V^+$ and ${\bold R}$ for $R_{gr}$. Note 
that if $\pi$ is a path with $l(\pi) \ge k$ and $v(\pi,k) = v_0v_1...v_{k-1}$ then
$|v_{k-1}| = |v_k| + 1 \ge 1.$  Thus $v(\pi,k) \in {\bold V}^k.$

\definition{Definition 4.1} $$Path_k =  span \{ v(\pi,k) \ | \ \pi \ is \ a \ path, 
l(\pi) \ge k \}$$  
and
$$Path_k(v) = v{\bold V}^{k-1} \cap Path_k .$$
\enddefinition
\bigskip
Let $f : {\bold V} \rightarrow F$ be defined by f(v) = 1 for all $v \in {\bold V}$ and  
 $I^l$ denote $I \otimes ... \otimes I$, taken $l$ times. Let $g_l :{\bold V}^l \rightarrow {\bold V}^{l-1}$ be 
defined by $g_l = f \otimes I^{l-1}$. 

For any vertex $v \in {\bold V}$ and any $l \geq 0$ define

$$S_l(v) = span \{ u \ | \ v > u, |u| = |v| - l \}$$

and

$$P_l(v) = span \{ u-w \ | \ v > u, v > w, |u| = |w| = |v| - l \}. $$
 
Note that $P_l(v) = S_l(v) = (0)$ if $l > |v|$, 
that $S_0(v) = span \{ v \}$, and that $P_0(v) = (0).$

Note also that $P_l(v) = ker f \mid_{S_l(v)}$ and therefore 
$$P_l(v){\bold V}^m = \ ker \ g_{m+1}\mid_{S_l(v){\bold V}^m}$$
for all $l,m \geq 0.$ Combining this with Proposition 3.6, we have 

$$gr \  {\bold R}_2 = span \{ v(u-w) \ | \ u,w \in {\Cal S}_1(v), v \in {\bold V} \}$$

$$= \sum_{v \in {\bold V}} vP_1(v) .$$  

\smallskip
\proclaim{Lemma 4.2} For $k \geq 2$, 
$${\bold R}^{(k)} = Path_k \ \cap \ \bigcap_{i=0}^{k-2} ker\ (I^{i+1} \otimes f \otimes I^{k-i-2}) .$$
\endproclaim

\smallskip
\demo{Proof} Since $${\bold R}^{(k)} = \bigcap_{i=0}^{k-2} {\bold V}^i{\bold R}{\bold V}^{k-i-2}$$
we have

$${\bold R}^{(k)\bot} = \sum_{i=0}^{k-2} {\bold V}^{*i}{\bold R}^\bot {\bold V}^{*k-i-2}.$$

We also have that 
$${\bold R}^\bot = span \{ \{ v^*u^*  \ | \ |u| \neq |v|-1 \ or \ v  \not> u \} \cup \{ v^* f \ | \ v \in {\bold V} \}\}$$

$$= span \{ v^*u^*  \ | \ |u| \neq |v|-1 \ or \ v \not> u \}  + {\bold V}^*f .$$

Let 
$$ M = span \{ v^*u^*  \ | \ |u| \neq |v|-1 \ or \ v \not> u \}.$$

Then 

$$ {\bold R}^{(k)\bot} = \sum_{i=0}^{k-2} \{{\bold V}^{*i}M{\bold V}^{*k-i-2} + {\bold V}^{*i+1}f{\bold V}^{*k-i-2} \}$$

$$ = \sum_{i=0}^{k-2} {\bold V}^{*i}M{\bold V}^{*k-i-2} + \sum_{i=0}^{k-2} {\bold V}^{*i+1}f{\bold V}^{*k-i-2}$$

$$ = ((\sum_{i=0}^{k-2} {\bold V}^{*i}M{\bold V}^{*k-i-2})^\bot \cap (\sum_{i=0}^{k-2} {\bold V}^{*i+1}f{\bold V}^{*k-i-2})^\bot)^\bot.$$
So
$${\bold R}^{(k)} = (\sum_{i=0}^{k-2} {\bold V}^{*i}M{\bold V}^{*k-i-2})^\bot  \cap (\sum_{i=0}^{k-2} {\bold V}^{*i+1}f{\bold V}^{*k-i-2})^\bot.$$

Now

$$(\sum_{i=0}^{k-2} {\bold V}^{*i}M{\bold V}^{*k-i-2})^\bot = \cap_{i=0}^{k-2} {\bold V}^iM^\bot {\bold V}^{k-i-2} = Path_k$$ 
and
$$(\sum_{i=0}^{k-2} {\bold V}^{*i}f{\bold V}^{*k-i-2})^\bot  =  \cap_{i=0}^{k-2} {\bold V}^{i+1}<f>^\bot {\bold V}^{k-i-2}$$

$$= \cap_{i=0}^{k-2} {\bold V}^{i+1}(ker f){\bold V}^{k-i-2} = \cap_{i=0}^{k-2}  ker (I^{i+1} \otimes f \otimes I^{k-i-2}),$$
giving the result.
\hfill $\square$ \enddemo

We will need the following result, whose proof is straightforward.

\smallskip
\proclaim{Lemma 4.3} Let $W_1$ and $W_2$ be $F$-vector spaces, $h :  W_1 \rightarrow W_2$  a linear transformation,
 $A \subseteq W_1$, $C \subseteq  W_2$, subspaces,  and $B = h^{-1}(C)$.  Then  
$$h(A) \cap  h(B) = h(A \cap B).$$
\endproclaim

\smallskip
\proclaim{Lemma 4.4} If $l \geq 0$, $j \geq 1$, and $v \in \cup_{j=2}^n  V_j$, then  
$$P_j(v){\bold V}^{l+1} \cap  {\bold R}^{(l+2)} = g_{l+3}(S_{j-1}(v){\bold V}^{l+2} \cap {\bold R}^{(l+3)}).$$
\endproclaim
\smallskip
\demo{Proof} Note that $$(f \otimes I)(S_{j-1}(v){\bold V} \cap  {\bold R}) \subseteq P_j(v)$$
  and   
$$(f \otimes I^{l+2})({\bold V}{\bold R}^{(l+2)}) \subseteq {\bold R}^{(l+2)}.$$
Consequently, 
$$g_{l+3}(S_{j-1}(v){\bold V}^{l+2} \cap {\bold R}^{(l+3)}) \subseteq P_j(v){\bold V}^{l+1} \cap  {\bold R}^{(l+2)}.$$
To prove the reversed inclusion we note that by Lemma 4.2, 
$$P_j(v){\bold V}^{l+1} \cap  {\bold R}^{(l+2)} = P_j(v){\bold V}^{l+1} \cap  Path_{l+2} \cap 
(\cap_{i=0}^l ker(I^{i+1} \otimes f \otimes I^{l-i}))$$ and
$$g_{l+3}(S_{j-1}(v){\bold V}^{l+2} \cap {\bold R}^{(l+3)}) = g_{l+3}(S_{j-1}(v){\bold V}^{l+2} \cap  Path_{l+3} \cap 
(\cap_{i=0}^{l+1} ker(I^{i+1} \otimes f \otimes I^{l+1-i}))).$$

Let $g = g_{l+4}...g_{l+j+2}$. Then for any subspace $W \subseteq {\bold V}^{l+2}$ we have
$$ g^{-1}(W) = {\bold V}^{j-1}W$$
and $$(g_{l+3}g)^{-1}(W) = {\bold V}^jW.$$

Then $$P_j(v){\bold V}^{l+1} \cap  Path_{l+2} \cap 
(\cap_{i=0}^l ker(I^{i+1} \otimes f \otimes I^{l-i}))=$$
$$(g_{l+3}g)(Path_{j+1}(v){\bold V}^{l+1}) \cap (g_{l+3}g)({\bold V}^j(ker \ ( f \otimes I^{l+1}))) \  \cap$$
$$(g_{l+3}g)({\bold V}^jPath_{l+2}) \cap(\cap_{i=0}^l {\bold V}^j(ker(I^{i+1}\otimes f \otimes I^{l-i}))) = $$
$$(g_{l+3}g)(Path_{j+1}(v){\bold V}^{l+1} \cap {\bold V}^jPath_{l+2} \cap (\cap _{i=0}^{l+1} {\bold V}^j(ker ( I^i \otimes f \otimes I^{l+1-i})))) = $$
$$(g_{l+3}g)(Path_{j+l+2}(v) \cap(\cap_{i=0}^{l+1} {\bold V}^j(ker(I^i \otimes f \otimes I^{l+1-i})))).$$

Similarly,
$$g_{l+3}(S_{j-1}(v){\bold V}^{l+2} \cap  Path_{l+3} \cap 
(\cap_{i=0}^{l+1} ker(I^{i+1} \otimes f \otimes I^{l+1-i}))) =$$  
$$g_{l+3}(g(Path_j(v){\bold V}^{l+2}) \cap g({\bold V}^{j-1}Path_{l+3}) \ \cap$$
$$(\cap_{i=0}^{l+1} g({\bold V}^{j-1}(ker(I^{i+1}\otimes f \otimes I^{l+1-i}))))) \  = $$
$$(g_{l+3}g)(Path_j(v){\bold V}^{l+2} \cap {\bold V}^{j-1}Path_{l+3} 
\cap(\cap_{i=0}^{l+1} {\bold V}^{j-1}(ker(I^{i+1} \otimes f \otimes I^{l+1-i})))) = $$
$$(g_{l+3}g)(Path_{j+l+2}(v) \cap(\cap_{i=0}^{l+1} {\bold V}^j(ker(I^i \otimes f \otimes I^{l+1-i})))),$$
proving the lemma.
\hfill $\square$ \enddemo

\medskip

\proclaim{Lemma 4.5} Suppose $\{P_1(v){\bold V}^k\} \cup \{{\bold V}^i{\bold R}{\bold V}^{k-i-1}| 0 \le i \le
k-1\}$ is distributive for any $v \in {\bold V}$. 
Then $\{{\bold V}^i{\bold R}{\bold V}^{k-i}|0 \le i \le k\}$ is distributive.
\endproclaim
\medskip
\demo{Proof} By Lemma 1.1 of \cite {SW} it is sufficient to
prove that 
$$\{v{\bold V}^{k+1} \cap {\bold V}^i{\bold R}{\bold V}^{k-i}|0 \le i \le k\}$$ 
is distributive for
all $v \in {\bold V}.$
Now $g_{k+2}$ restricts to  an isomorphism of $v{\bold V}^{k+1}$ onto ${\bold V}^{k+1}$.  Since
$g_{k+2}(v{\bold V}^{k+1} \cap {\bold R}{\bold V}^k) = P_1(v){\bold V}^k$ and
$g_{k+2}(v{\bold V}^{k+1} \cap {\bold V}^{i+1}{\bold R}{\bold V}^{k-i-1}) = {\bold V}^i{\bold R}{\bold V}^{k-i-1}$ for $0 \le i \le
k-1$, the result follows.
\hfill $\square$ \enddemo

\proclaim{Theorem 4.6} Let $\Gamma$ be a uniform layered graph
with a unique minimal element. Then 
$\{ {\bold V}^i{\bold R}{\bold V}^{k-1} \; | \; 0 \leq i \leq k \}$ 
generates a distributive  lattice in $T({\bold V})$. Consequently, 
$gr \ A(\Gamma)$ is a Koszul algebra.
\endproclaim
In view of Lemma 4.5, this will follow from:
\bigskip
\proclaim{Lemma 4.7} $\{ P_l(v){\bold V}^k\} \cup \{{\bold V}^i{\bold R}{\bold V}^{k-i-1} \; | \; 0 \le i \le
k-1 \}$ is  distributive  for all $k \ge 1$ and all $l > 0$.
\endproclaim

\demo{Proof}  The proof is by induction on $k$, the result being trivial
for $k = 1$.
We assume   $\{ P_l(v){\bold V}^m\} \cup \{{\bold V}^i{\bold R}{\bold V}^{m-i-1} \; | \; 0 \le i \le m-1 \}$
is  distributive  for all $m < k $ and all $l > 0$.

First note that any proper subset of $\{ P_l(v){\bold V}^k\} \cup \{{\bold V}^i{\bold R}{\bold V}^{k-i-1} \;
| \; 0 \le i \le k-1 \}$ is  distributive.
Indeed, by Lemma 4.5, $\{{\bold V}^i{\bold R}{\bold V}^{k-i-1} \; | \; 0 \le i \le k-1 \}$ is
distributive.  Hence it is sufficient to show that
$\{ P_l(v){\bold V}^k\} \cup \{{\bold V}^i{\bold R}{\bold V}^{k-i-1} \; | \; 0 \le i <j\} \cup
\{{\bold V}^i{\bold R}{\bold V}^{k-i-1} \; | \; j <  i \le k-1 \}$
is  distributive  for all $j, 0 \le j \le k-1\}.$ 
Now let ${\Cal K}_{j,1} = \{ P_l(v){\bold V}^j\} \cup \{{\bold V}^i{\bold R}{\bold V}^{j-i-1} \; | \; 0 \le
i \le j-1 \}$
and ${\Cal K}_{j,2} = \{ {\bold V}^i{\bold R}{\bold V}^{k-j-2-i} \; | \; 0 \le i \le k-j-2 \}$.
Then ${\Cal K}_{j,1}$ and ${\Cal K}_{j,2}$ are distributive
by the induction assumption.
Since  $\{ P_l(v){\bold V}^k\} \cup \{{\bold V}^i{\bold R}{\bold V}^{k-i-1} \; | \; 0 \le i <j\} \cup
\{{\bold V}^i{\bold R}{\bold V}^{k-i-1} \; | \; j <  i \le k-1 \} = {\Cal K}_{j,1}{\bold V}^{k-j} \cup
{\bold V}^{j+1}{\Cal K}_{j,2}$, the assertion follows.

In view of Theorem 1.2 of \cite{SW}, it is therefore sufficient to prove that
$$ (P_l(v){\bold V}^k \cap {\bold R}{\bold V}^{k-1} \cap ... \cap {\bold V}^i{\bold R}{\bold V}^{k-i-1})
\cap ({\bold V}^{i+1}{\bold R}{\bold V}^{k-i-2} + ... + {\bold V}^{k-1}{\bold R})$$
$$ = (P_l(v){\bold V}^k \cap {\bold R}{\bold V}^{k-1} \cap ... \cap {\bold V}^{i+1}{\bold R}{\bold V}^{k-i-2})
$$
$$+ (P_l(v){\bold V}^k \cap {\bold R}{\bold V}^{k-1} \cap ... \cap {\bold V}^i{\bold R}{\bold V}^{k-i-1} \cap
({\bold V}^{i+2}{\bold R}{\bold V}^{k-i-3} + ... + {\bold V}^{k-1}{\bold R})).$$

Now write 
$$X_i = S_l(v){\bold V}^k \cap {\bold R}{\bold V}^{k-1} \cap ... \cap {\bold V}^i{\bold R}{\bold V}^{k-i-1},$$ 
$$Y_i
= {\bold V}^{i+1}{\bold R}{\bold V}^{k-i-2},$$
and 
$$Z_i = {\bold V}^{i+2}{\bold R}{\bold V}^{k-i-3} + ... + {\bold V}^{k-1}{\bold R}.$$  
Then we need to show that

$$ ker \ g_{k+1}|_{X_i \cap (Y_i + Z_i)} = ker \ g_{k+1}|_{X_i \cap Y_i} +
ker \ g_{k+1}|_{X_i \cap Z_i}.$$

The right-hand side is contained in the left-hand side, so it is enough to
prove equality of dimensions.

Hence it is enough to prove

$$dim \ X_i \cap (Y_i + Z_i) - dim \ g_{k+1} (X_i \cap (Y_i + Z_i)) $$
$$ = dim \ X_i \cap Y_i - dim \ g_{k+1}(X_i \cap Y_i) + dim \ X_i \cap Z_i -
dim \ g_{k+1}(X_i \cap Z_i)$$
$$ - dim \ X_i \cap Y_i \cap Z_i + dim \ g_{k+1}(X_i \cap Y_i \cap Z_i).$$

Now, by the induction assumption and Lemma 4.5,  $\{X_i, Y_i, Z_i\}$ is
distributive.
Therefore, the desired equality is equivalent to
$$dim \ g_{k+1}(X_i \cap (Y_i + Z_i)) = dim \ g_{k+1}(X_i \cap Y_i) + $$
$$dim \ g_{k+1}(X_i \cap Z_i) - dim \ g_{k+1}(X_i \cap Y_i \cap Z_i).$$

But, since $\{X_i, Y_i, Z_i\}$ is distributive,
$$g_{k+1}(X_i \cap (Y_i + Z_i)) = g_{k+1}(X_i \cap Y_i + X_i \cap Z_i) =
g_{k+1}(X_i \cap Y_i) + g_{k+1}(X_i \cap Z_i).$$

Hence we need only show that $$dim \ g_{k+1}(X_i \cap Y_i) \cap g_{k+1}(X_i\cap Z_i)
 = dim \ g_{k+1}(X_i \cap Y_i \cap Z_i).$$
Since $$g_{k+1}(X_i \cap Y_i \cap Z_i) \subseteq 
g_{k+1}(X_i \cap Y_i) \cap g_{k+1}(X_i\cap Z_i)$$
this is equivalent to 
$$g_{k+1}(X_i \cap Y_i) \cap g_{k+1}(X_i \cap Z_i) = g_{k+1}(X_i \cap Y_i
\cap Z_i)$$

Now by Lemma 4.4, the left-hand side of this expression is equal to
$$P_{l+1}(v){\bold V}^{k-1}\cap {\bold R}{\bold V}^{k-2} \cap ... \cap {\bold V}^i{\bold R}{\bold V}^{k-i-2} \cap
({\bold V}^{i+1}{\bold R}{\bold V}^{k-i-3} + ... + {\bold V}^{k-2}{\bold R}).$$

In view of the distributivity of $\{S_l(v){\bold V}^k\} \cup \{{\bold V}^t{\bold R}{\bold V}^{k-t-1} | 0 \le
t \le k-1\}$,
which follows from Lemma 1.1 of \cite{SW} and the induction assumption,
the right-hand side of the expression may be written as
$$g_{k+1}(S_l(v){\bold V}^k \cap {\bold R}{\bold V}^{k-1} \cap ... \cap {\bold V}^{i+2}{\bold R}{\bold V}^{k-i-3}) $$
$$+ g_{k+1}(S_l(v){\bold V}^k \cap {\bold R}{\bold V}^{k-1} \cap ... \cap {\bold V}^{i+1}{\bold R}{\bold V}^{k-i-2} \cap
({\bold V}^{i+3}{\bold R}{\bold V}^{k-i-4} + ... + {\bold V}^{k-1}{\bold R})).$$

By Lemma 4.4, this is equal to
$$P_{l+1}(v){\bold V}^{k-1} \cap {\bold R}{\bold V}^{k-2} \cap ... \cap {\bold V}^{i+1}{\bold R}{\bold V}^{k-i-3} $$
$$+P_{l+1}(v){\bold V}^{k-1} \cap {\bold R}{\bold V}^{k-2} \cap ... \cap {\bold V}^i{\bold R}{\bold V}^{k-i-2} \cap
({\bold V}^{i+2}{\bold R}{\bold V}^{k-i-4} + ... + {\bold V}^{k-2}{\bold R}).$$

By the induction assumption, these expressions for the left-hand and right-hand 
sides are equal, so the proof is complete.
\hfill $\square$  \enddemo 

\bigskip

\head 5. Koszulity of $A(\Gamma)$ \endhead
\medskip

\proclaim {Theorem 5.1}Let $\Gamma$ be a uniform layered graph
with a unique minimal element. 
Then $A(\Gamma)$ is a Koszul algebra.
\endproclaim

\demo{Proof}
This follows from Theorem  4.6 and Proposition 7.2 of \cite{PP}.
\hfill $\square$ \enddemo


\Refs \widestnumber\key{GGRSW}

\ref\key GGR \by I. Gelfand, S. Gelfand, and V. Retakh
\paper Noncommutative algebras
associated to complexes and graphs  \jour
Selecta Math. (NS)  \vol 7 \yr 2001\pages 525--531  \endref

\ref\key GGRSW \by I. Gelfand, S. Gelfand, V. Retakh, S. Serconek,
and R. Wilson\paper Hilbert series of quadratic algebras
associated with decompositions of noncommutative polynomials \jour
J. Algebra \vol 254 \yr 2002\pages 279--299  \endref

\ref\key GRSW \by  I. Gelfand, V. Retakh, S. Serconek,  and R. Wilson
\paper On a class of algebras associated to directed graphs
\paperinfo math.QA/0506507  \jour Selecta Math. (N.S.) \vol 11
\yr 2005 \endref

\ref\key GRW \by I. Gelfand, V. Retakh, and R. Wilson \paper
Quadratic-linear algebras associated with decompositions of
noncommutative polynomials and differential polynomials \jour
Selecta Math. (N.S.) \vol 7\yr 2001 \pages 493--523 \endref

\ref\key Pi \by D. Piontkovski \paper  Algebras associated to pseudo-roots 
of noncommutative polynomials  are Koszul \paperinfo math.RA/0405375 \endref

\ref\key PP \by A. Polishchuk, and L. Positselski
\paper Quadratic Algebras
\paperinfo preprint
\yr 2005
\endref



\ref\key SW \by S. Serconek and R. L. Wilson \paper Quadratic algebras
associated with decompositions of noncommutative polynomials
are Koszul algebras \jour J. Algebra \vol 278 \yr 2004 \pages 473-493 \endref



\endRefs
\enddocument